\documentclass{article}

\usepackage{amssymb}
\usepackage[cp1251]{inputenc}
\usepackage[english]{babel}
\usepackage{amsmath}
\usepackage{color}
\usepackage{graphics}
\usepackage{epsfig}
\usepackage[all,knot]{xy}
\xyoption{arc}

\def\mapr#1{\smash{\mathop{\buildrel{#1}\over\longrightarrow}}}
\def\mapl#1{\smash{\mathop{\buildrel{#1}\over\longleftarrow}}}

\newtheorem{theorem}{Theorem}

\def\qed{\hfill\vrule width2mm height2mm depth2mm}

\def\proof{{\bf Proof.}}

\def\F{{\bf F}}

\def\R{{\bf R}}

\def\b{{\bf b}}

\def\e{{\bf e}}

\def\0{{\bf 0}}
\def\1{{\bf 1}}

\def\b#1{{\bf #1}}

\def\cA{{\cal A}}

\def\cG{{\cal G}}

\def\cK{{\cal K}}

\def\ad{{\hbox{\bf ad}}}

\def\bydef{\stackrel{def}{=}}

\def\Der{\mbox{\bf Der }}

\def\Hom{{\hbox{\bf Hom}\;}}

\def\Int{{\hbox{\bf Int}\;}}

\def\Mor{{\hbox{\bf Mor}\;}}
\def\Obj{\hbox{\bf Obj}\;}

\def\Out{{\hbox{\bf Out}\;}}

\def\fcolor{\color{black}}
\def\blf{\fcolor $}
\def\elf{$ \color{black}}

\def\bdf{\fcolor$$}
\def\edf{$$\color{black}}

\def\beq#1{\fcolor\begin{equation}\label{#1}}
\def\eeq{\end{equation}\color{black}}

\def\b${\color{blue} $}
\def\e${$ \color{black}\hskip-0.05cm}

\title{Derivations of Group Algebras and  Hochschild cohomology}

\author{A.~S.~Mishchenko\thanks{This work was supported by the Russian Foundation for Basic Research (Grant
No 18-01-00398)}}
\date{}
\begin{document}
\maketitle

\begin{abstract}
This article is a presentation of a report read at an international conference dedicated to the memory of Professor B.Yu. Sternin in Moscow on November 6-09, 2018.
The Hochschild homology and  cohomology group can be described in terms of the homology and cohomology of the classifying space of the groupoid of
the adjoint action of the group under the suitable assumption of the finiteness of the supports of cohomology groups.
The difference between homology and cohomology leads to a correction of the results in the book by D. J. Benson \cite{benson-1995},\cite{benson-1991}.
\end{abstract}

\section{Introduction}

This talk is motivated by comparing the results of our last work, Arutyunov A.A., Mishchenko A.S.,(2018), \cite{Ar-Mi-2018}, (the full text is presented in the journal "Matematicheskij sbornik" \cite{Ar-Mi-2018-2}) and also works by Arutyunov, A.A., Mishchenko, A.S., Shtern, A.I., (2016), \cite{Ar-Mi-St-2016}, \cite{Ar-Mi-St-2017},
in which a description of the algebra of exterior derivations of the group algebra \blf R[G]\elf of a finitely representable discrete group \blf G\elf is presented in terms of the Cayley complex of the groupoid \blf\cG\elf of the adjoint action of the group \blf G\elf, with the results of Burghelea (1985) \cite{burghelea-1985} and Benson (1995, 1991) \cite{benson-1995}, \cite{benson-1991}, which describe the Hochschild homology and cohomology of the group algebra \blf R[G]\elf in terms of the classifying spaces \blf BC\langle x\rangle\elf centralizers \blf C\langle x\rangle\elf of conjugate classes \blf \langle x\rangle\elf of groups \blf G\elf.

The space of external derivations has a description in the form of one-dimensional Hochschild cohomology of the same group algebra (see the book by R. Pierce (1986) \cite{Pierce-1986}, the definition "a", p. 248).
Therefore, a natural question arises: is it possible to describe all Hochschild cohomology of group algebra in terms of geometric constructions on the groupoid of the adjoint action of a group by analogy with the external derivations of group algebra?

In the book of Benson (1991, 1995) \cite{benson-1991}, \cite{benson-1995}
just contains the calculation of the Hochschild homology and cohomology of the group algebra \blf C[G]\elf in terms of centralizers of classes of conjugate elements in the group \blf G\elf. In particular, there is a theorem (Volume 2, p.76)

\begin{theorem}[2.11.2]
The additive structure of the Hochschild homology and cohomology of the group algebra $RG$ is given by the formulas
\begin{enumerate}
\item \blf HH_{n}(RG)\cong\bigoplus\limits_{g\in G^{G}} H_{n}(C_{G}(g),R)\elf,
\item \blf HH^{n}(RG)\cong\bigoplus\limits_{g\in G^{G}} H^{n}(C_{G}(g), R)\elf.
\end{enumerate}
\end{theorem}

This statement contradicts our statement (Arutyunov A.A., Mishchenko A.S.,(2018),\cite{Ar-Mi-2018-2},Corollary 1),
in which the one-dimensional Hochschild cohomology is described as the sum of one-dimensional cohomology with finite supports of the classifying centralizers of classes of conjugate elements of the group $G$, while Benson’s book deals on ordinary cohomology of the same spaces.

This is due to the fact that in the Benson’s book the definition of Hochschild cohomology is artificially modified and drawn to the duality of Hochschild homology. Namely, the Hochschild cohomology of the group algebra $RG$ Benson means not cohomology $HH^{n}(RG; RG)$, in which the algebra $RG$ itself serves as coefficients, as is customary in most of the works devoted to the Hochschild cohomology (see for example, the works \cite{felix-2004}, \cite{gerstenhaber-1963}, \cite{kontsevich-1999}). In order to reduce the entire description of the Hochschild cohomology to the Hochschild homology without relative to the use of Hochschild cohomology in specific problems, but only for the sake of duality of Hochschild homology and cohomology the dual bimodule $(RG)^{*}=\Hom(RG,R)$ is used. In other words, $HH^{n}(RG)$ denotes the group $HH^{n} (RG)=HH^{n}(RG,(RG)^{*})$, while under Hochschild homology
the standard definition of $HH_{n}(RG)=HH_{n}(RG,RG)$ is understood.

However, in Benson’s book, the proof of the Theorem 2.11.2 relies on dubious equality (on the same page):
\beq{e7}
R_{\Delta(G)}\Uparrow^{G\times G}\downarrow_{\Delta(G)}\cong
\prod\limits_{g\in G^{G}}R_{\Delta(C_{G}(g))}\Uparrow^{\Delta(G)},
\eeq
which, without any arguments, refers to Theorem 3.3.4 Mackey on the decomposition of arbitrary modules over group algebras (\cite{benson-1995}, p. 61).

The above formula (1) cannot be satisfied for finite (noncommutative groups) from dimensional considerations. Indeed, the left-hand side of the formula has the dimension
\bdf
\dim\left(
R_{\Delta(G)}\Uparrow^{G\times G}\downarrow_{\Delta(G)}
\right)=(\#G)^{2}.
\edf
Each right-hand factor has a dimension
\bdf
\dim\left(R_{\Delta(C_{G}(g))}\Uparrow^{\Delta(G)}\right)=\#G.
\edf
In total, the classes of conjugate elements for noncommutative groups are strictly less than the cardinality of the group.
\bdf
\#G^{G}<\#G.
\edf

This means that the dimension of the right-hand side is calculated as
\bdf
\dim\left( \prod\limits_{g\in G^{G}}R_{\Delta(C_{G}(g))}\Uparrow^{\Delta(G)}\right)=\#G^{G}\cdot\#G<
(\#G)^{2}.
\edf

So we have the inequality
\bdf
\dim\left(R_{\Delta(G)}\Uparrow^{G\times G}\downarrow_{\Delta(G)}\right)>
\dim\left(\prod\limits_{g\in G^{G}}R_{\Delta(C_{G}(g))}\Uparrow^{\Delta(G)}\right).
\edf

Note that in the original theorem of Burghelea (1985) \cite{burghelea-1985}
only half of the theorem from Benson’s book is stated. Namely, here is written:
\begin{theorem}[I, Burghelea, p.361]
  1) \blf HH_{*}(R[G])=\sum\limits_{\hat{x}\in\langle G\rangle}H_{*}(BG_x; R).\elf
\end{theorem}

At about the same time another book by Weibel (1997) \cite{weibel-1997}
was published, in which (Corollary 9.7.5) the Burghelea theorem on calculating the Hochschild homology of the group algebra in terms of centralizers of classes of conjugate elements was formulated. The calculation of the Hochschild cohomology, which are given in the book of Benson, is prudently omitted.

We suggest a uniform way of describing both the Hochschild homology and cohomology of the group algebra \blf C[G]\elf in terms of the classifying space \blf B\cG\elf of the groupoid \blf \cG\elf of the adjoint action of the group \blf G\elf. In these terms, the Hochschild homology of the group algebra \blf C[G]\elf coincides with the homology of the classifying space \blf B\cG \elf. The Hochschild cohomology can also be identified with the invariants of the classifying space \blf B\cG\elf of the groupoid \blf\cG\elf, namely, with the cohomology of this space, but with some finiteness conditions for cochains on \blf B\cG\elf.

\section{Hochschild Cohomology}
\subsection{Derivations}
Consider the Banach algebra \blf  \cA \elf  and a \blf  \cA \elf --bimodule
\blf  E \elf . A linear mapping \bdf D: \cA \mapr{}E \edf is called a derivation (or
differentiation) if, for any elements \blf  a,b\in\cA\elf , the so-called Leibniz
identity (with respect to the two-sided action of the algebra \blf \cA\elf   on the
bimodule \blf  E\elf ) \bdf D(ab)=D(a)b+aD(b), \quad a,b\in\cA. \edf holds (see
Definition~1.8.1 in the Dales paper (2000) \cite{Dales-2000}).

The space \blf \Der(\cA,E)\elf  of all derivations from \blf \cA\elf  to  \blf
E\elf  has the subspace of inner derivations \blf \Int(\cA,E)\subset\Der(\cA,E)\elf ,
defined by the adjoint representations \bdf \ad_{x}(a)\bydef xa-ax, \quad
x\in E, a\in\cA. \edf

The quotient space \blf \Out(\cA,E)= \Der(\cA,E)/\Int(\cA,E)\elf  is called  the space
of outer derivations; this space can be interpreted using the one-dimensional
Hochschild cohomology of the algebra \blf \cA\elf  with coefficients in the bimodule \blf  E\elf :
\bdf HH^{1}(\cA;E)\approx\Out(\cA,E), \edf

We have proved
that the algebra of outer derivations
\blf \Out(C[G])=\Der(C[G])/\Int(C[G])\elf  of the algebra \blf C[G]\elf  is isomorphic to the
one-dimensional cohomology of the Cayley complex \blf \cK(\cG)\elf  of the groupoid \blf \cG\elf
with finite supports:
\bdf HH^{1}(C[G])\approx\Out(C[G])\approx H^{1}_{f}(\cK(\cG); \R). \edf

\subsection{Hochschild cohomology}
So we want to express the Hochschild cohomology \blf HH^{k}(C[G])\elf in the terms of geometric properties of the groupoid \blf\cG\elf.\vskip 1cm

Let \blf\Lambda\elf be a group algebra of the group \blf G\elf, \blf\Lambda =C[G]\elf.
The Hochschild cochain complex consists of \blf C^{k}(\Lambda)\elf spaces of multilinear mappings
\bdf
f:\underbrace{\Lambda\times\Lambda\times\cdots\times\Lambda}_{k}\mapr{}\Lambda.
\edf

Any such mapping is completely given by the formula:
\bdf
f(g_{1},g_{2},\dots,g_{k})=\sum\limits_{h\in G}f^{h}_{g_{1},g_{2},\dots,g_{k}}h,
\quad f^{h}_{g_{1},g_{2},\dots,g_{k}}\in \mathbb{R},
\edf
 for \blf g_{1},g_{2},\dots,g_{k}\in G\elf, moreover, the matrix satisfies the condition \blf\left\|f^{h}_{g_{1},g_{2},\dots,g_{k}}\right\|\elf satisfies the finiteness condition:
 \begin{itemize}
\item[(\F)] For any set of items \blf g_{1},g_{2},\dots,g_{k}\in G\elf the set of indexes \blf h\in G\elf, for which \blf x^{h}_{g}\neq 0\elf is finite:
\bdf
\{h\in G: f^{h}_{g_{1},g_{2},\dots,g_{k}}\neq 0\}<+\infty.
\edf
\end{itemize}

In fact, if
\blf x=(x^{1},x^{2},\dots,x^{k})\in\underbrace{\Lambda\times\Lambda\times\cdots\times\Lambda}_{k}\elf is arbitrary element, \blf x^{j}\in\Lambda\elf,  and
\blf
x^j=\sum\limits_{g_{j}\in G}\lambda^{j}_{g_{j}}\cdot g_{j}, \label{e2}
\elf
then
\bdf
\begin{array}{l}
f(x)=f(x^{1},x^{2},\dots,x^{k})=\\\\=f\left(\left(\sum\limits_{g_{1}\in G}\lambda^{1}_{g_{1}}\cdot g_{1}\right),\left( \sum\limits_{g_{2}\in G}\lambda^{2}_{g_{2}}\cdot g_{2}\right),\cdots,
\left(\sum\limits_{g_{k}\in G}\lambda^{k}_{g_{k}}\cdot g_{k}\right)\right)=\\\\=
\sum\limits_{g_{1},g_{2},\dots,g_{k}\in G}\lambda^{1}_{g_{1}}\lambda^{2}_{g_{2}}\cdots
\lambda^{k}_{g_{k}}f(g_{1},g_{2},\dots,g_{k})=\\=
\sum\limits_{g_{1},g_{2},\dots,g_{k},h\in G}\lambda^{1}_{g_{1}}\lambda^{2}_{g_{2}}\cdots
\lambda^{k}_{g_{k}}f^{h}_{g_{1},g_{2},\dots,g_{k}} h.
\end{array}
\edf
The finiteness condition \blf(\F)\elf of the matrix \blf\left\|f^{h}_{g_{1}, g_{2}, \dots, g_{k}}\right \|\elf guarantees the finite sums in the formulas .

The Hochschild cochain complex
\bdf
0\mapr{}C^{0}(\Lambda)\mapr{\partial_{0}}C^{1}(\Lambda)\mapr{\partial_{1}}C^{2}(\Lambda)\mapr{\partial_{2}}
\cdots\mapr{\partial_{k-1}}C^{k}(\Lambda)\mapr{\partial_{k}}C^{k+1}(\Lambda)\mapr{\partial_{k+1}}\cdots
\edf
is defined by the formula for \blf f\in C^{k}(\Lambda)\elf, \blf \partial_{k}(f)\in C^{k+1}(\Lambda)\elf:
\bdf
\begin{array}{l}
\partial_{k}(f)(g_{1},g_{2},\dots,g_{k+1})=\\\\
=g_{1}f(g_{2},\dots,g_{k+1})-f(g_{1}g_{2},g_{3},\dots,g_{k+1})+
f(g_{1},g_{2}g_{3},g_{4}\dots,g_{k+1})+\\
+\cdots+(-1)^{j}f(g_{1},g_{2},\dots,g_{j}g_{j+1}, g_{j+2},\dots,g_{k+1})+\cdots+\\+(-1)^{k}f(g_{1},g_{2},\dots,g_{k}g_{k+1})
+(-1)^{k+1}f(g_{1},g_{2},\dots,g_{k})g_{k+1}=\\\\
=g_{1}f(g_{2},\dots,g_{k+1})+\\
+\sum\limits_{j=1}^{k}(-1)^{j}f(g_{1},\dots, g_{j-1},g_{j}g_{j+1},g_{j+2},\dots,g_{k+1})+\\
+(-1)^{k+1}f(g_{1},g_{2},\dots,g_{k})g_{k+1}.
\end{array}
\edf

Every $k$ -dimensional simplex $\sigma$ of the classifying space $B\cG$ of the groupoid $\cG$ is a sequence of morphisms
\bdf
\sigma:\left(a_{0}\mapr{g_{k}}a_{1}\mapr{g_{k-1}}\cdots\mapr{g_{k-j+1}}a_{j}
\mapr{g_{k-j}}a_{j+1}\mapr{g_{k-j-1}}\cdots\mapr{g_{1}}a_{k}\right),
\edf
and
\bdf
a_{j+1}=g_{k-j}a_{j}g_{k-j}^{-1},
\edf
in particular
\bdf
a_{k}=(g_{1}g_{2}\dots g_{k})a_{0}(g_{1}g_{2}\dots g_{k})^{-1}.
\edf

Well, then $k$ - dimensional cochains $T\in C^{k}(B\cG)$ on $B\cG$ are such functions
$T$ on simplices $\sigma$ that are given by their values $T(\sigma)$:
\bdf
T(\sigma)=T_{g_{1},g_{2},\dots,g_{k}}^{a_{0}}.
\edf
It turns out cochain complex:
\bdf
0\mapr{}C^{0}(B\cG)\mapr{\delta_{0}}C^{1}(B\cG)\mapr{\delta_{1}}\cdots
\mapr{\delta_{k-1}}C^{k}(B\cG)\mapr{\delta_{k}}C^{k+1}(B\cG)\mapr{\delta_{k+1}}\cdots
\edf
Now we will construct a mapping of cochain complexes:
\begin{equation}\label{2}
\xymatrix{
0\ar[r]&C^{0}(C[G])\ar[r]^{\partial_{0}}\ar[d]_{T_{0}}&
C^{1}(C[G])\ar[r]^{\partial_{1}}\ar[d]_{T_{1}}&
\cdots\\
0\ar[r]&C^{0}(B\cG)\ar[r]^{\delta_{0}}&C^{1}(B\cG)\ar[r]^{\delta_{1}}&
\cdots\\
\cdots\ar[r]^{\partial_{k-1}}&C^{k}(C[G])\ar[r]^{\partial_{k}}\ar[d]_{T_{k}}&
C^{k+1}(C[G])\ar[r]^{\partial_{k+1}}\ar[d]_{T_{k+1}}&\cdots\\
\cdots\ar[r]^{\delta_{k-1}}&C^{k}(B\cG)\ar[r]^{\delta_{k}}&
C^{k+1}(B\cG)\ar[r]^{\delta_{k+1}}&\cdots
}
\end{equation}
The maps  $T_{k}$ are constructed by the formula: if $f\in C^{k}(\cA)$,
\bdf
f:\underbrace{C[G]\times C[G]\times\cdots\times C[G]}_{k}\mapr{} C[G],
\edf

\bdf
f(g_{1},g_{2},\dots,g_{k})=\sum\limits_{h\in G}f^{h}_{g_{1},g_{2},\dots,g_{k}}h,
\edf
and the simplex \b$\sigma\e$ has the form
\bdf
\sigma=\left(a_{0}\mapr{g_{k}}a_{1}\mapr{g_{k-1}}\cdots\mapr{g_{k-j+1}}a_{j}
\mapr{g_{k-j}}a_{j+1}\mapr{g_{k-j-1}}\cdots\mapr{g_{1}}a_{k}\right),
\edf
such that
\bdf
a_{j+1}=g_{k-j}a_{j}g_{k-j}^{-1},
\edf
then
\bdf
T_{k}(f)(\sigma)=f^{h}_{g_{1},g_{2},\dots,g_{k}},
\edf
where
\bdf
h=
(g_{1}g_{2}\dots g_{k})a_{0}=a_{k}(g_{1}g_{2}\dots g_{k}).
\edf

\begin{theorem}
The diagram (\ref{2}) above is commutative and induces the isomorphism in cohomology under some condition of finiteness.
\bdf
T_k:HH^k(C[G],C[G])\mapr{}H^k_f(B\cG).
\edf
\end{theorem}
\proof

Let $f\in C^{k}(\Lambda)$ be an arbitrary element. In according with definition the element $f$ is a polylinear function

$$
f:\underbrace{\Lambda\otimes\Lambda\otimes\cdots\otimes\Lambda}_{k}\mapr{}\Lambda,
$$
that is defined by the formula for arbitrary arguments $u_{1},u_{2},\dots,u_{k}\in\Lambda$,
$$
u_{j}=\sum\limits_{g_{j}\in G} \lambda^{g_{j}}_{j}\cdot g_{j}:
$$
$$
\begin{array}{l}
f(u_{1}\otimes u_{2}\otimes\dots\otimes u_{k})=\\=f\left(\sum\limits_{g_{1}\in G} \lambda^{g_{1}}_{1}\cdot g_{1}\otimes
\sum\limits_{g_{2}\in G} \lambda^{g_{2}}_{2}\cdot g_{2}\otimes\cdots\otimes
\sum\limits_{g_{k}\in G} \lambda^{g_{k}}_{k}\cdot g_{k}\right)=\\=
\sum\limits_{(g_{1}\otimes g_{2}\otimes\dots\otimes g_{k})}\lambda_{1}^{g_{1}}
\lambda_{2}^{g_{2}}\cdots\lambda_{k}^{g_{k}}
f(g_{1}\otimes g_{2}\otimes\dots\otimes g_{k})
\end{array}
$$

We need to establish the identity
$$
T_{k+1}\partial_{k}(f)\equiv \delta_{k}T_{k}(f)
$$
for any function
$$
f=f(g_{1},g_{2},\dots,g_{k})=\sum\limits_{h\in G}f^{h}_{g_{1}g_{2}\dots g_{k}}\cdot h.
$$
Identity checking needs on every simplex
$$
\sigma=\left(a_{0}\mapr{g_{k+1}}a_{1}\mapr{g_{k}}\cdots\mapr{g_{k-j+1}}a_{j}
\mapr{g_{k-j+1}}a_{j+1}\mapr{g_{k-j}}\cdots\mapr{g_{1}}a_{k+1}\right).
$$

\subsubsection{The left hand part:}
$$
T_{k+1}\partial_{k}(f)(\sigma)=(\partial_{k}(f))^{h}_{g_{1},g_{2},\dots,g_{k+1}},
$$
where
$$
h=(g_{1}g_{2}\cdots g_{k+1})a_{0}=a_{k+1}(g_{1}g_{2}\cdots g_{k+1}).
$$
Further,
$$
\partial_{k}(f)(g_{1},g_{2},\dots,g_{k+1})=\sum\limits_{h\in G} (\partial_{k}(f))^{h}_{g_{1},g_{2},\dots,g_{k+1}}\cdot h,
$$
$$
\begin{array}{l}
\partial_{k}(f)(g_{1},g_{2},\dots,g_{k+1})=\\=
g_{1}f(g_{2},\dots,g_{k+1})-f(g_{1}g_{2},\dots,g_{k+1})+\\+
f(g_{1},g_{2}g_{3},\dots,g_{k+1})-f(g_{1},g_{2},g_{3}g_{4},\dots,g_{k+1})+\cdots\\+
(-1)^{k}f(g_{1},g_{2},g_{3}g_{4},\dots,g_{k}g_{k+1})+
(-1)^{k+1}f(g_{1},g_{2},g_{3}g_{4},\dots,g_{k})g_{k+1}=\\\\
=\sum\limits_{h\in G}\left(
g_{1}hf^{h}_{g_{2},\dots,g_{k+1}}-f^{h}(g_{1}g_{2},\dots,g_{k+1})h+\right.\\+
f^{h}(g_{1},g_{2}g_{3},\dots,g_{k+1})h-f^{h}(g_{1},g_{2},g_{3}g_{4},\dots,g_{k+1})h+
\cdots\\+\left.
(-1)^{k}f^{h}(g_{1},g_{2},g_{3}g_{4},\dots,g_{k}g_{k+1})h+
(-1)^{k+1}f^{h}(g_{1},g_{2},g_{3}g_{4},\dots,g_{k})hg_{k+1}\right)=\\\\=
\sum\limits_{h\in G}\left(
hf^{g_{1}^{-1}h}_{g_{2},\dots,g_{k+1}}-f^{h}_{g_{1}g_{2},\dots,g_{k+1}}h+\right.\\+
f^{h}_{g_{1},g_{2}g_{3},\dots,g_{k+1}}h-f^{h}_{g_{1},g_{2},g_{3}g_{4},\dots,g_{k+1}}h+
\cdots\\+\left.
(-1)^{k}f^{h}_{g_{1},g_{2},g_{3}g_{4},\dots,g_{k}g_{k+1}}h+
(-1)^{k+1}f^{hg_{k+1}^{-1}}_{g_{1},g_{2},g_{3}g_{4},\dots,g_{k}}h\right)=\\\\=
\sum\limits_{h\in G}\left(
f^{g_{1}^{-1}h}_{g_{2},\dots,g_{k+1}}-f^{h}_{g_{1}g_{2},\dots,g_{k+1}}+\right.\\+
f^{h}_{g_{1},g_{2}g_{3},\dots,g_{k+1}}-f^{h}_{g_{1},g_{2},g_{3}g_{4},\dots,g_{k+1}}+
\cdots\\+\left.
(-1)^{k}f^{h}_{g_{1},g_{2},g_{3}g_{4},\dots,g_{k}g_{k+1}}+
(-1)^{k+1}f^{hg_{k+1}^{-1}}_{g_{1},g_{2},g_{3}g_{4},\dots,g_{k}}\right) h=\\\\=
\sum\limits_{h\in G}(\partial_{k}(f))^{h}_{g_{1},g_{2},\dots,g_{k+1}} h.
\end{array}
$$
Hence
$$
\begin{array}{l}
T_{k+1}\partial_{k}(f)(\sigma)=(\partial_{k}(f))^{h}_{g_{1},g_{2},\dots,g_{k+1}}=\\\\=
f^{g_{1}^{-1}h}_{g_{2},\dots,g_{k+1}}-f^{h}_{g_{1}g_{2},\dots,g_{k+1}}+\\+
f^{h}_{g_{1},g_{2}g_{3},\dots,g_{k+1}}-f^{h}_{g_{1},g_{2},g_{3}g_{4},\dots,g_{k+1}}+
\cdots\\+
(-1)^{k}f^{h}_{g_{1},g_{2},g_{3}g_{4},\dots,g_{k}g_{k+1}}+
(-1)^{k+1}f^{hg_{k+1}^{-1}}_{g_{1},g_{2},g_{3},g_{4},\dots,g_{k}},
\end{array}
$$
where
$$
h=(g_{1}g_{2}\cdots g_{k+1})a_{0}=a_{k+1}(g_{1}g_{2}\cdots g_{k+1}).
$$

\subsubsection{Right hand part:}

$$
\begin{array}{l}
\delta_{k}T_{k}(f)(\sigma)=\sum\limits_{j=0}^{k+1}(-1)^{j}T_{k}(f)(d_{j}\sigma)=\\\\=
f^{h'}_{g_{1},g_{2},\dots,g_{k}}+
\sum\limits_{j=1}^{k}(-1)^{j}f^{h}_{g_{1},\dots,g_{k-j+1}g_{k-j+2},\dots,g_{k},g_{k+1}}+\\\\+
(-1)^{k+1}f^{h''}_{g_{2},g_{3},\dots,g_{k+1}}
\end{array}
$$
where
$$
\begin{array}{l}
h=(g_{1}g_{2}\dots g_{k+1})a_{0}=a_{k+1}(g_{1}g_{2}\dots g_{k+1}),\\
hg_{k+1}^{-1}=h'=(g_{1}g_{2}\dots g_{k})a_{1}=a_{k+1}(g_{2}\dots g_{k+1}),\\
g_{1}^{-1}h=h''=(g_{2}\dots g_{k+1})a_{0}=a_{k}(g_{2}\dots g_{k+1})
\end{array}
$$

This means that

$$
T_{k+1}\partial_{k}(f)\equiv (-1)^{k+1}\delta_{k}T_{k}(f)
$$
\qed

\subsection{Hochschild Homology}
In contrast to the cohomology, the Hochschild homology is constructed using the Hochschild chain complex (see definition Benson (1991) \cite{benson-1991}, p. 74)
\bdf
0\mapl{}C_{0}(\Lambda)\mapl{d_{1}}C_{1}(\Lambda)\mapl{d_{2}}
C_{2}(\Lambda)\mapl{d_{3}}\cdots\mapl{d_{k}}
C_{k}(\Lambda)\mapl{d_{k+1}}C_{k+1}(\Lambda)\mapl{d_{k+2}}\cdots
\edf

Here
\bdf
\begin{array}{l}
C_{0}(\Lambda)=\Lambda,\\
C_{1}(\Lambda)=\Lambda\otimes\Lambda,\\
C_{2}(\Lambda)=\Lambda\otimes\Lambda\otimes\Lambda,\\
\vdots\\
C_{k}(\Lambda)=\underbrace{\Lambda\otimes\Lambda\otimes\cdots\otimes\Lambda}_{k+1}.\\
\vdots
\end{array}
\edf

And

\bdf
\begin{array}{lll}
d_{1}(g_{0}\otimes g_{1})&=&g_{0}g_{1}-g_{1}g_{0},\\
&&\quad g_{0}\otimes g_{1}\in C_{1}(\Lambda),\\
d_{2}(g_{0}\otimes g_{1}\otimes g_{2})&=& g_{0}g_{1}\otimes g_{2}-
g_{0}\otimes g_{1}g_{2}+g_{2}g_{0}\otimes g_{1},\\
&&\quad g_{0}\otimes g_{1}\otimes g_{2}\in C_{2}(\Lambda),\\
\vdots\\
d_{k}(g_{0}\otimes g_{1}\otimes\cdots\otimes g_{k})&=&
g_{0}g_{1}\otimes\cdots\otimes g_{k}-\\&&-
g_{0}\otimes g_{1}g_{2}\otimes\cdots\otimes g_{k}+\cdots+\\&&+
(-1)^{k-1}g_{0}\otimes g_{1}\otimes g_{2}\otimes\cdots\otimes g_{k-1}g_{k}+\\
&&+(-1)^{k}g_{k}g_{0}\otimes g_{1}\otimes\cdots\otimes g_{k-1},\\
&&\quad g_{0}\otimes g_{1}\otimes\cdots\otimes g_{k}\in C_{k}(\Lambda),\\
\vdots
\end{array}
\edf

\subsubsection{The chain complex of classifying space of the groupoid $\cG$}
The groupoid chain complex $\cG$
\bdf
0\mapl{}C_{0}(B\cG)\mapl{\delta_{1}}C_{1}(B\cG)\mapl{\delta_{2}}\cdots
\mapl{\delta_{k}}C_{k}(B\cG)\mapl{\delta_{k+1}}C_{k+1}(B\cG)\mapl{\delta_{k+2}}\cdots
\edf
is induced by linear combinations of simplices
\bdf
\sigma=\left(a_{0}\mapr{g_{k+1}}a_{1}\mapr{g_{k}}\cdots\mapr{g_{k-j+2}}a_{j}
\mapr{g_{k-j+1}}a_{j+1}\mapr{g_{k-j}}\cdots\mapr{g_{1}}a_{k+1}\right)\in C_{k+1}(B\cG)
\edf such that
$$
g_{k-j+1}a_{j}=a_{j+1}g_{k-j+1},
$$
where
$$
\begin{array}{l}
\delta_{k+1}(\sigma)=\delta_{k+1}\left(a_{0}\mapr{g_{k+1}}a_{1}\mapr{g_{k}}\cdots\mapr{g_{k-j+2}}a_{j}
\mapr{g_{k-j+1}}a_{j+1}\mapr{g_{k-j}}\cdots\mapr{g_{1}}a_{k+1}\right)=\\\\
=\left(a_{1}\mapr{g_{k}}\cdots\mapr{g_{k-j+2}}a_{j}
\mapr{g_{k-j+1}}a_{j+1}\mapr{g_{k-j}}\cdots\mapr{g_{1}}a_{k+1}\right)-\\
-\left(a_{0}\mapr{g_{k}g_{k+1}}a_{2}\mapr{g_{k-1}}\cdots\mapr{g_{k-j+2}}a_{j}
\mapr{g_{k-j+1}}a_{j+1}\mapr{g_{k-j}}\cdots\mapr{g_{1}}a_{k+1}\right)+\cdots+\\
+(-1)^{j}\left(a_{0}\mapr{g_{k+1}}a_{1}\mapr{g_{k}}\cdots\mapr{g_{k-j+3}}a_{j-1}
\mapr{g_{k-j+1}g_{k-j+2}}a_{j+1}\mapr{g_{k-j}}\cdots\mapr{g_{1}}a_{k+1}\right)+\cdots+\\
+(-1)^{k+1}\left(a_{0}\mapr{g_{k+1}}a_{1}\mapr{g_{k}}\cdots\mapr{g_{k-j+2}}a_{j}
\mapr{g_{k-j+1}}a_{j+1}\mapr{g_{k-j}}\cdots\mapr{g_{2}}a_{k}\right)
\end{array}.
$$
By analogy with cohomology, we construct a map of chain complexes
\begin{equation}\label{3}
\xymatrix{
0&C_{0}(\Lambda)\ar[l]\ar[d]_{S_{0}}&
C_{1}(\Lambda)\ar[l]_{d_{1}}\ar[d]_{S_{1}}&
\cdots\ar[l]_{d_{2}}
\\
0\ar[r]&C_{0}(B\cG)&C_{1}(B\cG)\ar[l]^{\delta_{1}}&
\cdots\ar[l]^{\delta_{2}}\\
\cdots&C_{k}(\Lambda)\ar[l]_{d_{k}}\ar[d]_{S_{k}}&
C_{k+1}(\Lambda)\ar[l]_{d_{k+1}}\ar[d]_{S_{k+1}}&\cdots\ar[l]_{\hbox{\phantom{-----}}d_{k+2}}
\\
\cdots&C_{k}(B\cG)\ar[l]^{\delta_{k}}&
C_{k+1}(B\cG)\ar[l]^{\delta_{k+1}}&\cdots\ar[l]^{\hbox{\phantom{-----}}\delta_{k+2}}
}
\end{equation}

Let \blf g_{0}\otimes g_{1}\otimes g_{2}\otimes\cdots\otimes g_{k}\in C_{k}(\Lambda)\elf.
Put
\bdf
\begin{array}{l}
S_{k}(g_{0}\otimes g_{1}\otimes g_{2}\otimes\cdots\otimes g_{k})=\\\\=
\left(a_{0}\mapr{g_{k}}a_{1}\mapr{g_{k-1}}\cdots\mapr{g_{k-j+1}}a_{j}
\mapr{g_{k-j}}a_{j+1}\mapr{g_{k-j-1}}\cdots\mapr{g_{1}}a_{k}\right),
\end{array}\edf
where
\bdf
a_{0}=g_{0}g_{1}\dots g_{k}.
\edf
In particular, if \blf g_{0}\in C_{0}(\Lambda)\elf, then \blf S_{0}(g_{0})=g_{0}\in C_{0}(B\cG)\elf.
If
\blf g_{0}\otimes g_{1}\in C_{1}(\Lambda)\elf, then

\bdf
S_{1}(g_{0}\otimes g_{1})=\left(a_{0}\mapr{g_{1}}a_{1}\right)\in C_{1}(B\cG), \quad a_{0}=g_{0}g_{1},\quad a_{1}=g_{1}g_{0}.
\edf

\begin{theorem}
The diagram (\ref{3}) above is commutative and generates the isomorphism in homologies:
\bdf
S_k:HH_k(C[G])\mapr{}H_k(B\cG).
\edf
\end{theorem}
\proof  \   Let $g_{0}\otimes g_{1}\otimes\cdots\otimes g_{k+1}\in C_{k+1}(\Lambda)$. We should check
$$
S_{k}d_{k+1}(g_{0}\otimes g_{1}\otimes\cdots\otimes g_{k+1})=
\delta_{k+1}S_{k+1}(g_{0}\otimes g_{1}\otimes\cdots\otimes g_{k+1}).
$$
\subsubsection{The left hand part:}

$$
\begin{array}{lll}
S_{k}d_{k+1}(g_{0}\otimes g_{1}\otimes\cdots\otimes g_{k+1})&=&
S_{k}\left(g_{0}g_{1}\otimes\cdots\otimes g_{k+1}\right)-\\
&&-
S_{k}\left(g_{0}\otimes g_{1}g_{2}\otimes\cdots\otimes g_{k+1}\right)+\cdots+\\
&&+
(-1)^{k}S_{k}\left(g_{0}\otimes g_{1}\otimes g_{2}\otimes\cdots\otimes g_{k}g_{k+1}\right)+\\
&&+S_{k}\left((-1)^{k+1}g_{k+1}g_{0}\otimes g_{1}\otimes\cdots\otimes g_{k}\right),
\end{array}
$$
or
$$
\begin{array}{l}
S_{k}d_{k+1}(g_{0}\otimes g_{1}\otimes\cdots\otimes g_{k+1})=\\\\=
S_{k}\left(g_{0}g_{1}\otimes\cdots\otimes g_{k+1}\right)-\\-
S_{k}\left(g_{0}\otimes g_{1}g_{2}\otimes\cdots\otimes g_{k+1}\right)+\cdots+\\
+
(-1)^{k}S_{k}\left(g_{0}\otimes g_{1}\otimes g_{2}\otimes\cdots\otimes g_{k}g_{k+1}\right)+\\
+(-1)^{k+1}S_{k}\left(g_{k+1}g_{0}\otimes g_{1}\otimes\cdots\otimes g_{k}\right)=\\\\=
\left(a_{0}\mapr{g_{k+1}}a_{1}\mapr{g_{k}}\cdots\mapr{g_{k-j+2}}a_{j}
\mapr{g_{k-j+1}}a_{j+1}\mapr{g_{k-j}}\cdots\mapr{g_{2}}a_{k}\right)-\\
\phantom{aaaaaaaaa}\left\{\hbox{Here }a_{0}=(g_{0}g_{1}g_{2}g_{3}\cdots g_{k+1})\right\}\\\\
-\left(a_{0}\mapr{g_{k+1}}a_{1}\mapr{g_{k}}\cdots\mapr{g_{k-j+2}}a_{j}
\mapr{g_{k-j+1}}a_{j+1}\mapr{g_{k-j}}\cdots\mapr{g_{1}g_{2}}a_{k}\right)+\\
\phantom{aaaaaaaaa}\left\{\hbox{Here }a_{0}=(g_{0}g_{1}g_{2}g_{3}\cdots g_{k+1})\right\}\\
+\dots+
\\\\
+(-1)^{k}\left(a_{0}\mapr{g_{k}g_{k+1}}a_{1}\mapr{g_{k-1}}\cdots\mapr{g_{k-j+1}}a_{j}
\mapr{g_{k-j}}a_{j+1}\mapr{g_{k-j-1}}\cdots\mapr{g_{1}}a_{k}\right)+\\
\phantom{aaaaaaaaa}\left\{\hbox{Here }a_{0}=(g_{0}g_{1}g_{2}g_{3}\cdots g_{k+1})\right\}\\\\
+(-1)^{k+1}\left(a_{0}\mapr{g_{k}}a_{1}\mapr{g_{k-1}}\cdots\mapr{g_{k-j+1}}a_{j}
\mapr{g_{k-j}}a_{j+1}\mapr{g_{k-j-1}}\cdots\mapr{g_{1}}a_{k}\right).\\
\phantom{aaaaaaaaa}\left\{\hbox{Here }a_{0}=g_{k+1}g_{0}(g_{1}g_{2}g_{3}\cdots g_{k})\right\}.
\end{array}
$$

$$
\begin{array}{l}
S_{k}(g_{0}\otimes g_{1}\otimes g_{2}\otimes\cdots\otimes g_{k})=\\\\=
\left(a_{0}\mapr{g_{k}}a_{1}\mapr{g_{k-1}}\cdots\mapr{g_{k-j+1}}a_{j}
\mapr{g_{k-j}}a_{j+1}\mapr{g_{k-j-1}}\cdots\mapr{g_{1}}a_{k}\right),
\end{array}
$$

\bdf
a_{0}=(g_{0}g_{1}g_{2}\dots g_{k}).
\edf

\subsubsection{The right hand part:}

$$
\begin{array}{l}
\delta_{k+1}S_{k+1}(g_{0}\otimes g_{1}\otimes\cdots\otimes g_{k+1})=\\\\
=\delta_{k+1}\left(a_{0}\mapr{g_{k+1}}a_{1}\mapr{g_{k}}\cdots\mapr{g_{k-j+2}}a_{j}
\mapr{g_{k-j+1}}a_{j+1}\mapr{g_{k-j}}\cdots\mapr{g_{1}}a_{k+1}\right)
\end{array}
$$
where

\bdf
a_{0}=(g_{0}g_{1}g_{2}\dots g_{k+1}).
\edf
Hence

$$
\begin{array}{l}
\delta_{k+1}S_{k+1}(g_{0}\otimes g_{1}\otimes\cdots\otimes g_{k+1})=\\=
\delta_{k+1}\left(a_{0}\mapr{g_{k+1}}a_{1}\mapr{g_{k}}\cdots\mapr{g_{k-j+2}}a_{j}
\mapr{g_{k-j+1}}a_{j+1}\mapr{g_{k-j}}\cdots\mapr{g_{1}}a_{k+1}\right)=\\\\
=\left(a_{1}\mapr{g_{k}}\cdots\mapr{g_{k-j+2}}a_{j}
\mapr{g_{k-j+1}}a_{j+1}\mapr{g_{k-j}}\cdots\mapr{g_{1}}a_{k+1}\right)-\\
-\left(a_{0}\mapr{g_{k}g_{k+1}}a_{2}\mapr{g_{k-1}}\cdots\mapr{g_{k-j+2}}a_{j}
\mapr{g_{k-j+1}}a_{j+1}\mapr{g_{k-j}}\cdots\mapr{g_{1}}a_{k+1}\right)+\cdots+\\
+(-1)^{j}\left(a_{0}\mapr{g_{k+1}}a_{1}\mapr{g_{k}}\cdots\mapr{g_{k-j+3}}a_{j-1}
\mapr{g_{k-j+1}g_{k-j+2}}a_{j+1}\mapr{g_{k-j}}\cdots\mapr{g_{1}}a_{k+1}\right)+\cdots+\\
+(-1)^{k+1}\left(a_{0}\mapr{g_{k+1}}a_{1}\mapr{g_{k}}\cdots\mapr{g_{k-j+2}}a_{j}
\mapr{g_{k-j+1}}a_{j+1}\mapr{g_{k-j}}\cdots\mapr{g_{2}}a_{k}\right)
\end{array}.
$$

The left hand side coincides with the right hand side.

\qed
\section{Classifying space $B\cG$ of the groupoid $\cG$.}

\subsection{Right action: Classifying space $Br\cG$ of groupoid $r\cG$}

 The groupoid \blf r\cG\elf is defined as
\blf \Obj(r\cG)\approx G\elf, \blf\Mor(a,b)=\{x:ax=b\}, \quad a,b\in G\elf.

\begin{theorem}
Classifying space of groupoid \blf Br\cG\elf is contractible:
\bdf
\pi_k(Br\cG)=0, k\geq 0.
\edf
\end{theorem}

\subsection{Trivial action: Classifying space $BG$ of group $G$}

The group \blf G\elf is the category with single object and
\blf\Mor(G)\approx G\elf.
\begin{theorem}
Classifying space of group \blf BG\elf is the Eilenberg-Maclane complex:
\bdf
BG\sim K(G,1):
\edf
\bdf
BG \approx (Br\cG)/G.
\edf
\end{theorem}

\subsection{Adjoint action:Classifying space of groupoid $B\cG$}
\begin{theorem}
Classifying space of groupoid \blf B\cG\elf is didjoint union
\bdf
B\cG\approx \coprod\limits_{\langle g\rangle\in\langle G\rangle}B\cG_{\langle g\rangle}
\edf
\bdf
B\cG_{\langle g\rangle}\approx (Br\cG)/C(\langle g\rangle)\sim K(C(\langle g\rangle),1)
\edf
\end{theorem}

\subsection*{Acknowledgments}

When writing this paper, Shtern Alexander Isakovich provided invaluable assistance, to whom we express our sincere thanks. Thanks also to Lazarev Andrei Borisovich for valuable advice.

This work was supported by the Russian Foundation for Basic Research (Grant
No 18-01-00398)


\begin{thebibliography}{ccc}

\bibitem{benson-1995}
    {\sf  D. J. Benson,}
    {\it Representations and Cohomology, I. Basic representation theory of finite groups and associative algebras.}
    {\rm Cambridge University Press, 1995}

\bibitem{benson-1991}
    {\sf  D. J. Benson,}
    {\it  Representations and Cohomology, II. Cohomology of groups and modules.}
    {\rm Cambridge University Press, 1991}



\bibitem{Ar-Mi-2018}
    {\sf Arutyunov A.A., Mishchenko A.S.,}
    {\it Smooth Version of Johnson's Problem Concerning Derivations of Group Algebras}
    {\rm arXive e-print, t. 1801, No. 03480, p. 1-52 (2018)}


\bibitem{Ar-Mi-2018-2}
 Arutyunov,~A.\,A.,  Mishchenko,~A.\,S.
Smooth Version of Johnson's Problem Concerning Derivations of Group Algebras.
\textit{arXiv:1801.03480 [math.AT]}, (Submitted to Mathematical sbornik).


\bibitem{Ar-Mi-St-2016}
    {\sf Arutyunov A.A., Mishchenko A.S., Stern A.I.,}
    {\it Derivations of Group Algebras}
    {\rm Fundamental and applied mathematics, t. 21, issue 6, 2016, p. 63--75}



\bibitem{Ar-Mi-St-2017}
    {\sf Arutyunov A.A., Mishchenko A.S., Stern A.I.,}
    {\it Derivations of Group Algebras}
    {\rm arXive e-print, t. 1708, No. 05005, p. 1-28 (2017)}



\bibitem{burghelea-1985}
    {\sf Dan Burghelea}
    {\it  The cyclic homology of the group rings,}
    {\it  Comment. Math. Helvetici 60 (1985) 354-365 }

\bibitem{Pierce-1986}
    {\sf R.~S.~Pierce,}
    {\it Associative Algebras}, {\rm Springer-Verlag, New York--Berlin, 1982. --- 448 pp.}


\bibitem{gerstenhaber-1963}
    {\sf M. Gerstenhaber}
    {\it  The Cohomology Structure of an Associative Ring,}
    {\it  The Annals of Mathematics,  Second Series, Vol. 78, No. 2(Sep., 1963), pp. 267-288}

\bibitem{felix-2004}
    {\sf Y.Felix, J.-C. Thjmas, M. Vigue-Poirrier,}
    {\it The Hochschild cohomology of a closed manifold,}
    {\rm Publications Mathematiques de l'Institut des Hautes Etudes Scientifiques. Volume 99, Issue 1, 2004, pp 235–252.}

\bibitem{kontsevich-1999}
    {\sf M.Kontsevich,}
    {\it Operads and Motives in Deformation Quantization,}
    {\rm Letters in Mathematical Physics 48: 35-72, 1999. }



\bibitem{weibel-1997}
    {\sf  C.A.Weibel,}
    {\it  An Introduction to Homological Algebra.}
    {\rm Cambridge University Press, 1997}






\end{thebibliography}
\end{document}